\documentclass[12pt]{article}

\author{Moshchevitin, N.G.
\footnote{ Research is supported by the grants RFFI
(02--01--00192), grant of the President of RF (MD-3321.2004.1) and
INTAS (03-51-5070)}}
\title{The best Diophantine approximations: the phenomenon
of degenerate dimension}
\date{}
\begin{document}
\maketitle

This brief survey deals with  multi-dimensional Diophantine approximations in  sense of
linear form and with simultaneous Diophantine approximations. We discuss the phenomenon of
degenerate dimension of linear subspaces generated by the best Diophantine approximations.
Originally most of these results have been  established by the author in
\cite{MREGCH,MLAG,MDAN,BORD,MOSHEROT}. Here we collect  all of them together and give some
new formulations. In  contrast to our previous survey \cite{BORD}, this paper contains a
wider number of results, especially dealing with the best Diophantine approximations. It
also includes proofs or sometimes the  sketches of proofs. Some applications of these
results and  methods to the theory of small denominators  can be found  in
\cite{MREGCH,MOSHEGRAZ} and \cite{LOCHACK}.

\vskip+0.3cm

{\bf \S 1.\,\,\, The best Diophantine approximations in  sense of linear form.}

\vskip+0.3cm

{\bf 1.1 \hspace{2mm}  Notation.} \vskip+0.3cm

\par Let
$ \alpha _1,\dots ,\alpha _r  $ be real numbers which, together with 1, are linearly
independent over rationals. For an integer point
$$
m = ( m_0, m_1, \dots , m_r ) \in  {\bf Z}^{r+1} \setminus \{ (0,\dots ,0)\}
$$
we define
$$
\zeta (m) = | m_0 + m_1\alpha _1 + \dots + m_r \alpha _r | \hspace{2mm} \hspace{2mm} {\rm
and } \hspace{2mm} \hspace{2mm} M = \max_{j = 0, 1, \dots , r } |m_j|.
$$

A point
 $ m \in {\bf Z}^{r+1}\setminus \{ 0 \} $
is defined to be { \it the best approximation (in sense of linear form) } if
$$
\zeta (m) = \min_ { n \in Z^{r+1}\setminus \{ 0 \}: \hspace{2mm}  N \le M } \hspace{2mm} |
\zeta  (n) |
$$
(here $ N = \max_j |n_j|$). For the set of all best approximations $ m $ the corresponding
values of $ \zeta (m) $ and $ M $ can be ordered in descending (ascending) order:
$$
\zeta _1 > \zeta _2 > \dots > \zeta _\nu  > \zeta _{ \nu +1} > \dots ,
$$
$$
M_1 < M_2 < \dots < M_\nu  < M_{\nu +1} < \dots.
$$
(Here $ m_\nu = ( m_{0,\nu} , \dots , m_{r,\nu })   $ is
 $\nu -$th best approximation
and $ \zeta _\nu  = \zeta (m_\nu ),$ $  M_\nu  = \max_j |m_{j,\nu }|.$) By
 the
Minkowski convex body theorem it follows that $ \zeta _\nu  M_{\nu +1}^r \le 1 $. Let
 $ \Delta _\nu^r   $
denote the determinant of the
 $ r +1 $
consecutive best approximations:
$$
 \Delta_\nu^r  =
 \left|
\begin{array} {cccc}
m_{0,\nu } & m_{1,\nu } &  \dots & m_{r,\nu } \cr \dots & \dots & \dots & \dots \cr
m_{0,\nu +r} & m_{1,\nu +r} & \dots & m_{r,\nu +r}
\end{array}
  \right|.
$$
\vskip+0.3cm

{\bf 1.2. \hspace{2mm} The results on dimension.} \vskip+0.3cm

Here we  observe some properties of the  values $ \Delta_\nu^r $ discovered in \cite{MDAN}.

The following statement is well known from the continued fractions theory (see
\cite{Hin}). \vskip+0.3cm

\par {\bf Theorem 1.1. }
\hspace{2mm} {\it Let $ r = 1  $ and let $\alpha _1$ be an irrational number. Then for any
natural $ \nu  $ the determinant
 $ \Delta_\nu^1
$ is equal to
 $(-1)^{\nu -1}.
$} \vskip+0.3cm
\par

 The next result
deals with  dimension  $2$. It follows from the Minkowski convex body theorem.

\par {\bf Theorem 1.2. }
\hspace{2mm} {\it Let $ r = 2  $ and let $ \alpha _1,\alpha _2  $ be together with $1$
 linearly independent
over rationals. Then there exist  infinitely many values of $ \nu $
 for which} $ \Delta_\nu^2 \neq 0.
$ \vskip+0.3cm

As  was mentioned Theorem 1.2 is a simple corollary of the Minkowski theorem and we shall
give a sketch of the proof of Theorem 1.2 in the next section.

Now we    formulate our main result in this area which deals with the case $ r > 2$.

\vskip+0.3cm

\par {\bf Theorem 1.3. }
\hspace{2mm} {\it Given $ r \ge 3 $ there exists an uncountable set of $r$-tuples $
(\alpha _1,\dots ,\alpha _r )   $  such that  the corresponding sequence of the best
approximations
  $ m_\nu $
for all large  $ \nu  $ lies in a three-dimensional  sublattice $ \Lambda  (\alpha
_1,...,\alpha _r) $ of the lattice $ {\bf Z}^{r+1}.$
 Moreover,  each of these $r$-tuples
consists of   linearly independent over rationals together with $1$ reals. }

\vskip+0.3cm
\par {\bf Corollary.} \hspace{2mm}
For any $r$-tuple $ (\alpha _1,\dots ,\alpha _r )  $ in Theorem 1.3 there exists $\nu
_0(\alpha )$ such that for all $\nu >\nu _0(\alpha )$ we have
 $ \Delta_\nu^r = 0 .
$ \vskip+1.0cm

We shall give  the proof of the Theorem 1.3  in  Section 1.5. It is based on the
 so-called
Hinchin's singular $r$-tuples
 (see \cite{HINS1},\cite{HINS}
and Cassels' book \cite{CASS}). Before this proof in Section 1.4 we discuss the properties
of Hinchin's singular systems and their generalizations.

To finish this section we would like to emphasize once again that for $r\ge 2 $ for any
$r$-tuple $ (\alpha _1,\dots ,\alpha _r ) $ of $\bf Q$-independent reals all but finite
number  of the best approximation vectors never lie in a two-dimensional subspace but can
lie in a tree-dimensional subspace. We also would like to mention that there are many
results related to various definitions, algorithmic calculating  the best approximations
in general and for the algebraic numbers (see for example \cite{ Cus1, Cus2, Lag1, Lag2,
Lag3, Vor}) . \vskip+0.5cm

{\bf 1.3.\hspace{2mm}
 Sketch of the proof of Theorem 1.2.}
\vskip+0.3cm

\par
Assume the contrary: suppose that for some
 $ \alpha _1,\alpha _2 , $ which together with
$1$
 are linearly
independent over ${\bf  Z} $  we know that all best approximations
 $ m_\nu , \hspace{2mm}  \nu  \ge \nu _0 $
lie in some two-dimensional linear subspace $ \pi  $. Then from the
 continued fractions
theory (compare with Theorem 1.1) we have
\begin{equation}
 \zeta _\nu M_{\nu +1} \asymp
 \zeta _\mu M_{\mu +1} , \hspace{2mm}  \forall
 \nu ,\mu > \nu _0.
        \label{n1}
\end{equation}
Consider the cube
 $ E^3_H =\{ x = (x_1,x_2,x_3 ) \in R^3 :\hspace{2mm}  |x_j| \le H \} $
and the domain $ L^2_\sigma = \{ x \in R^3 ; \hspace{2mm}  \rho ( x; L^2) \le \sigma \}.$
Now the intersection
 $ \Omega( \sigma ,H) = E^3_H \cap L^2_\sigma  $
is a convex $O$-symmetric body  in $\pi $. As $ m_\nu  $ is a best approximation we
conclude that in the set
 $ \Omega ( \zeta _\nu ,M_{\nu +1} ) $
there are no integer points.
However from (\ref{n1}) it follows that
$$
Vol \hspace{2mm}
 \Omega ( \zeta _\nu ,M_{\nu +1} ) \asymp
\zeta _\nu  M^2_{\nu +1} \asymp M_{\nu +1} \to \infty,\hspace{2mm} \nu  \to \infty,
$$
and we have arrived at  a contradiction with the Minkowski convex body theorem.

\vskip+0.3cm

{\bf 1.4. \hspace{2mm}  Hinchin's $\psi$-singular linear forms.}

\par From the continued fractions theory \cite{Hin}
we know that in the case
 $ r = 1 $
we have
\begin{equation}
           \zeta _\nu  M_{\nu +1} \asymp 1.
           \label{n2}
\end{equation}
Next we show that  for linear forms in two or more variables
 the situation may be different:  the values of
$M_{\nu +1}$, corresponding to the best approximations $m_{\nu +1}$ may not be estimated
from below in terms of
 previous approximation $\zeta _\nu $. \vskip+0.3cm
\par{\bf Theorem 1.4 .}
\hspace{2mm} {\it Let $ r \ge 2 $ and  $ \psi (y ) $ be a real valued function decreasing
to
 $ 0 $  as
        $ y \to \infty. $
Then there exists an uncountable set of vectors $ (\alpha _1,\dots ,\alpha _r )  $ (with
components $ \alpha _1,\dots ,\alpha _r , 1  $ linearly independent over rationals) such
that for all corresponding best approximations we have  }
\begin{equation}
 \zeta _\nu  \le \psi ( M_{\nu +r- 1})\hspace{2mm} \forall\, \nu .
 \label{n3xaxa}
\end{equation}
\vskip+0.3cm We would like to remind  the readers the definition of a $\psi $-singular
linear form (in Hinchin's sense) \cite{HINS}, \cite{HINS1}. Let $\psi (y) = o(y^{-r}), y
\to + \infty$ decreases to zero. An $r$-tuple $ \alpha _1,\dots ,\alpha _r  $ is
$\psi$-singular   (in  sense of linear form) if for any $T> 1$ the following Diophantine
inequality has a solution in integer $r$-tuple $m$:
$$
|| m_1 \alpha _1+\dots + m_r\alpha _r || <\psi (T),\hspace{2mm} 0< \max_{1\le j \le r}
|m_j| \le T   .
$$
(Here $||\cdot ||$ denotes the distance to the nearest integer.)

 It is easy to verify that
an
 $r$-tuple $ \alpha _1,\dots ,\alpha _r $ is $\psi $-singular
 if and only if for all natural $\nu $
\begin{equation}
 \zeta _\nu  \le \psi ( M_{\nu +1}).
 \label{n3}
\end{equation}
From this point of view Theorem 1.4 in the case $ s\ge 3$ establishes the
 existence of
$r$-tuples which are "more singular" than Hinchin's singular linear forms.

\vskip+0.3cm

\par Proof of  Theorem  1.4.
\vskip+0.3cm

\par This proof was sketched in \cite{MDAN}.
It is based on the following lemma.

Let
$$
\nu_* \equiv \nu \pmod{r},\hspace{2mm}
1\le \nu _*\le r
$$
and
$\sigma $
is a gret positive number
($\sigma $
depends on
 $r$,
all the constants in symbols
$O(\cdot ), \ll ,\gg $
below may depend on
 $\sigma $).
Let
$$
\sigma _{j,\nu}=
      \sigma \nu _*^{j}, \hspace{2mm}
W =\max_{j,\nu }\sigma_{j,\nu}.
$$

\par{\bf Lemma 1.1.}  \hspace{2mm}
 {\it
There exists
an uncountable
set of vectors
$ \alpha _1,\dots ,\alpha _r  $
such that
\par 1) \hspace{2mm}
$1, \alpha _1,\dots ,\alpha _r  $
are linarly independent over $ {\bf Z} $;
\par 2) \hspace{2mm}
there exist a sequence of naturals
 $ p_\nu  $ under the conditions
\par (i) \hspace{2mm} $
        \sigma _{j,\nu}
       \psi (p_\nu ) \le
       || p_\nu  \alpha _j || =
       p_\nu  \alpha _j - a_{j,\nu }  \le
      (\sigma_{j,\nu }+1)
       \psi (p_\nu ) ;
       \hspace{2mm}  j = 1, \dots , r ,$
\par (ii) \hspace{2mm} $
      p_{\nu +1} \asymp p_\nu  (\psi(p_\nu ))^{-1}  .$}

\vskip+0.3cm

\par Proof of Lemma 1.1.

\vskip+0.3cm
\par We  construct numbers $\alpha_1 ,\dots ,\alpha_r$
with simultaneous approximations of  special type.

Let us take a sequence of zeroes and ones $\lambda=\{\lambda_\nu^2,\dots , \lambda_\nu^s
\}_{\nu=1}^\infty$; $\lambda_\nu^j\in\{0;1\}$. Now we define naturals $p_\nu$ ¨
$a_{j,\nu}$,
 $j=1,\dots , r$,
and segments $\Delta_{1,\nu}, \dots ,\Delta_{r,\nu}$ with lengths
$|\Delta_{j,\nu}|=2\psi(p_\nu)/p_\nu$ by the following recursive procedure.

Numbers $p_0$, $a_{1,0}, \dots , a_{r,0}$ may be taken arbitrary. Define
$$
\Delta_{j,0}=\biggl[\frac{a_{j,0}}{p_0}+
 \sigma _{j,\nu }
\frac{\psi(p_0)}{p_0},
\frac{a_{j,0}}{p_0}+
( \sigma _{j,\nu }  +1)
\frac{\psi(p_0)}{p_0}\biggr], \qquad
j=1,...,r.
$$
Let
$p_0,\dots,p_\nu$; $a_{j,0},\dots,a_{j,\nu}$ and
$\Delta_{j,0},\dots, \Delta_{j,\nu}$ are already defined.
We are constructing $p_{\nu+1}$, $a_{j,\nu+1}$ and $\Delta_{j,\nu+1}$.

Let
$$
p_{\nu+1}=\Bigl[6 p_\nu\bigl(\psi(p_\nu)\bigr)^{-1}\Bigr]+1.
$$
THen in any interval
of the length
 $\psi(p_\nu)/(6p_\nu)$
one can find a number
 $a/p_{\nu+1}$, $a\in {\bf Z}$.
Let
$$
\frac{a_{j,\nu+1}^0}{p_{\nu+1}}  \in
\biggl[\frac{a_{j,\nu}}{p_\nu}
+
( \sigma _{j,\nu }   +\frac16  )
\cdot\frac{\psi(p_\nu)}{p_\nu},
\frac{a_{j,\nu}}{p_\nu}
 +
 (\sigma _{j,\nu }
 +\frac26 )
  \cdot\frac{\psi(p_\nu)}{p_\nu}\biggr],
$$
$$
\frac{a_{j,\nu+1}^1}{p_{\nu+1}}  \in
\biggl[\frac{a_{j,\nu}}{p_\nu}+
 (\sigma _{j,\nu }
+
\frac46 ) \cdot\frac{\psi(p_\nu)}{p_\nu},
\frac{a_{j,\nu}}{p_\nu}+
 (\sigma _{j,\nu }
+
\frac56 )
\cdot\frac{\psi(p_\nu)}{p_\nu}\biggr].
$$
Now define
$$
\Delta_{j,\nu+1}^\tau =\biggl[\frac{a_{j,\nu+1}^\tau }{p_{\nu+1}}
+
\sigma _{j,\nu+1}
\frac{\psi(p_{\nu+1})}{p_{\nu+1}},
\frac{a_{j,\nu+1}^\tau }{p_{\nu+1}}+
(\sigma _{j,\nu +1}  +1)
\frac{\psi(p_{\nu+1})}{p_{\nu+1}}\biggr],
$$
$$
 \tau =0,1; j = 2, \dots , r.
$$
WE must note that
$$
\Delta_{j,\nu+1}^0\cap\Delta_{2,\nu+1}^1= \emptyset
$$
and
$$
\frac{a_{j,\nu}}{p_\nu}\notin\Delta_{j,\nu+1}^\tau , \qquad
\tau =0,1.
$$
Moreover $|\Delta_{j,\nu+1}^j|=\psi(p_{\nu+1})/p_{\nu+1}$
and due to
$$
(\sigma _{j,\nu +1}+1)
\psi(p_{\nu+1})/p_{\nu+1}\le\psi(p_\nu)/6p_\nu
$$
(here we suppose  $\psi$
to decrese fast enough: $W\psi (p_{\nu +1}) < \psi (p_\nu )$),
one has
$\Delta_{j,\nu+1}^\tau \subset\Delta_{j,\nu}$, $\tau =0,1$.
Put
$a_{j,\nu+1}=a_{j,\nu+1}^{\lambda_\nu^j}$
and $\Delta_{j,\nu+1}=\Delta_{j,\nu+1}^{\lambda_\nu}$.
An integer number
$a_{1,\nu+1}$
we define from the condition
$$
\frac{a_{1,\nu+1}}{p_{\nu+1}}\in
\biggl[\frac{a_{j,\nu}}{p_\nu}
+
( \sigma _{j,\nu}   +\frac16  )
\cdot\frac{\psi(p_\nu)}{p_\nu},
\frac{a_{j,\nu}}{p_\nu}
 +
 (\sigma _{j,\nu}
 +\frac26 )
  \cdot\frac{\psi(p_\nu)}{p_\nu}\biggr],
$$
Now let
$$
\Delta_{1,\nu+1}=\biggl[\frac{a_{1,\nu+1}}{p_{\nu+1}}
 \sigma _{j,\nu+1}
\frac{\psi(p_{\nu+1})}{p_{\nu+1}},
\frac{a_{1,\nu+1}}{p_{\nu+1}}
+
( \sigma _{j,\nu+1}   +1)
\frac{\psi(p_{\nu+1})}{p_{\nu+1}}\biggr].
$$
We have
 $\Delta_{1,\nu+1}\subset\Delta_{1,\nu}$,
$|\Delta_{1,\nu+1}|=\psi(p_{\nu+1})/p_{\nu+1}$ and
$
\frac{a_{1,\nu}}{p_\nu}\notin\Delta_{1,\nu+1}.
$

To summsrize we constructed a sequence of
enclosed segments
$$\{\Delta_{1,\nu}\}_{\nu=0}^\infty$$
and for arbitrary 0,1-sequence
 $\lambda$
we have   a sequence of enclosed segments
$\{\Delta_{j,\nu}\}_{\nu=0}^\infty$. Denote
$$
\alpha_1=\bigcap_\nu\Delta_{1,\nu}; \qquad
\alpha_j=\alpha_j(\lambda)=\bigcap_\nu\Delta_{j,\nu}.
$$
IN the sequences of fractions
$a_{1,\nu}/p_\nu$ and
$a_{j,\nu}/p_\nu$
all elements are different so
  $\alpha_1$, $
\alpha_j(\lambda)\notin {\bf Q}$.
Moreover we can choose $\lambda$ in such a way that
$1$, $\alpha_1 , \dots ,\alpha_r $
be linearly independent over rationals.
Similar procedure was
performed in \cite{ROD1qw}.

So we construct $\alpha_1,\dots ,\alpha_r\in {\bf R}$ satisfying the conditions

\vskip+0.2cm
\par {\bf A.} \,\,\, $1$, $\alpha_1$,..., $\alpha_r$ are linear independent ${\bf Z}$;

\par {\bf B.} \,\,\, for a sequence of naturals $p_\nu$,
\par $\|p_\nu\alpha_j\|=\|p_\nu\alpha_j-a_{j,\nu}\|<\psi(p_\nu)$, \
$j=1,\dots , r$,
\par
$3p_\nu\bigl(\psi(p_\nu)\bigr)^{-1}\le p_{\nu+1} \le4p_\nu\bigl(\psi(p_\nu)\bigr)^{-1}$ .
\vskip+0.2cm

One can easily verify that for any decreasing
 $\psi (y)$ the set
$$
M_\psi = \bigl\{(\alpha_1, \dots ,\alpha_r)\in {\bf R}^r:\alpha_1,\dots ,\alpha_r \,\,\,
{\rm satisfy}\,\,\, {\it  A, B}\bigr\}
$$
is uncountable and dense in  ${\bf R}^r$.

We must note for the future that the determinant
$$
 \left|
\begin{array} {cccc}
\sigma _{1,\nu }
 &
\sigma _{2,\nu }
 &  \dots &
\sigma _{r,\nu }
 \cr
\dots  & \dots & \dots &\dots \cr
\sigma _{1,\nu +r-1}
&
\sigma _{2,\nu +r-1}
 &  \dots &
\sigma _{r,\nu +r-1}
\end{array}
  \right| =
\pm \sigma ^r\prod_{1\le u<v\le r }(v-u)
$$
and for large
 $\sigma  = \sigma (r) $
for any
$\eta_{j,\mu} \in [-1,1],
j = 1,...,r; \mu = \nu ,...,\nu +r-1$
one has
$$
 \left|
\begin{array} {cccc}
\sigma _{1,\nu } +\eta_{1,\nu }
 &
\sigma _{2,\nu } +\eta_{2,\nu }
 &  \dots &
\sigma _{r,\nu }  +\eta_{r,\nu }
 \cr
\dots  & \dots & \dots &\dots \cr
\sigma _{1,\nu +r-1}  +\eta_{1,\nu +r-1}
&
\sigma _{2,\nu +r-1}
 +\eta_{2,\nu +r-1}
 &  \dots &
\sigma _{r,\nu +r-1}
 +\eta_{r,\nu +r-1}
\end{array}
  \right| =
$$
$$
=
\pm \sigma ^r
+ o(\sigma ^r) \neq 0.
$$
I the sequel this large  value of
$\sigma  = \sigma (r)$
is fixed.

Moreover it is easy to modify the construction in lemma 1.2
to establish that for any
value $p_\nu $ there exists
the {\it best } approximation for linear form
$||m^* _1 \alpha _1+...+m_r^* \alpha _r$
where  vectors
$(m_1^*,...,m^*_r)$,
$(p_\nu , a_{1,\nu },0 ,... , 0)$,
...,
$(p_\nu , 0 ,... , 0, \alpha _{r,\nu })$
are linearly dependent and $M^*\ll p^2_\nu $.

\par The proof is complete.
\vskip+0.3cm

\par{\bf Lemma 1.2.} \hspace{2mm}
{\it For numbers $\alpha_1 , \dots ,\alpha_r$, constructed in Lemma 1.1 there exists an
infinite sequence of values of the linear form
$$ \zeta (n_\nu)=n_{0,\nu}+n_{1,\nu}\alpha_1+\dots +n_{r,\nu}\alpha_r
=|n_{1,\nu}\alpha_1+\dots +n_{r,\nu}\alpha_r| $$ such that
$$
0<\zeta (n_{\nu+1})<\zeta (n_\nu)\ll\psi (\kappa N_{\nu+s-1}), $$ where
 $$
N_\nu=\max |n_{j,\nu}|  $$ and $ \kappa >0 $ is a constant.}

\vskip+0.3cm
\par Proof of Lemma 1.2.
\vskip+0.3cm

\par Let
$$
\zeta ( n_\nu  )
= \pm
 \left|
\begin{array} {cccc}
1 & \alpha _1 &  \dots & \alpha _r \cr
p_\nu  & a_{1,\nu } &  \dots & a_{r,\nu } \cr
\dots & \dots & \dots & \dots \cr
p_{\nu +r-1}  & a_{1,\nu +r-1 } &  \dots & a_{r,\nu +r-1}
\end{array}
  \right| =
$$
$$
= \pm
 \left|
\begin{array} {ccc}
  a_{1,\nu } -p_\nu \alpha _1 &  \dots & a_{r,\nu } -p_\nu  \alpha _r \cr
\dots & \dots & \dots \cr
 a_{1,\nu +r-1 } -p_{\nu +r-1}\alpha _1
&  \dots & a_{r,\nu +r-1}      -p_{\nu +r-1}\alpha _r
\end{array}
  \right| =
$$
$$
=
\psi (p_\nu )
 \cdots
\psi (p_{\nu+r-1} )
\times
$$
$$
\times
\left|
\begin{array} {cccc}
\sigma _{1,\nu } +\eta_{1,\nu }
 &
\sigma _{2,\nu } +\eta_{2,\nu }
 &  \dots &
\sigma _{r,\nu }  +\eta_{r,\nu }
 \cr
\dots  & \dots & \dots &\dots \cr
\sigma _{1,\nu +r-1}  +\eta_{1,\nu +r-1}
&
\sigma _{2,\nu +r-1}
 +\eta_{2,\nu +r-1}
 &  \dots &
\sigma _{r,\nu +r-1}
 +\eta_{r,\nu +r-1}
\end{array}
  \right|
\asymp
$$
$$
\asymp
\psi (p_\nu )
 \cdots
\psi (p_{\nu+r-1} ).
$$
(Here
$
\eta_{ j,\mu} =
\frac{a_{j,\mu}-p_\mu\alpha _j}{\psi (p_\mu) } +
\sigma _{j,\mu} \in [-1,1] $ and the sign
$+$ or $-$ is taken to satisfy $\zeta (n_\nu)>0$.)
Then
\begin{equation}
\zeta ( n_\nu)\asymp
\prod_{\mu = \nu}^{\nu +s-1}\max_{j=1,r}|a_{j,\mu}-p_{\mu}\alpha_j|
\asymp
\prod_{\mu = \nu}^{\nu +r-1}
\psi(p_\mu) <\psi(p_{\nu+r-1}).
\label{v1}
\end{equation}
For the coefficients $n_{j,\nu}$  we have
$$
n_{j,\nu}  =\pm
\left|
 \begin{array}{ccccccc}
p_\nu &
 a_{1,\nu}
     &
   \dots
     &
 a_{j-1,\nu}
     &
 a_{j+1,\nu}
     &
    \dots
     &
 a_{r,\nu}
        \cr
\dots &
 \dots
     &
   \dots
     &
 \dots
     &
 \dots
     &
    \dots
     &
 \dots
        \cr
p_\nu &
 a_{1,\nu +r-1}
     &
   \dots
     &
 a_{j-1,\nu +r-1}
     &
 a_{j+1,\nu +r-1}
     &
    \dots
     &
 a_{s,\nu +r-1}
\end{array}
\right| =
$$
$$
=\pm
\left|
\begin{array}{cccc}
p_\nu
   & a_{1,\nu}-p_\nu\alpha_1
   & \dots
   & a_{r,\nu}-p_\nu\alpha_r
   \cr \dots&\dots&\dots&\dots \cr
   p_{\nu  +r-1 }
   & a_{1,\nu +r-1}-p_{\nu +r-1}\alpha_1
   & \dots
   & a_{r,\nu +r-1}-p_{\nu+r-1}\alpha_r
\end{array}
\right|.
$$
Using (i) and (ii) we deduce
\begin{equation}
|n_{j,\nu +r-1}|\ll p_{\nu +r-1 }\psi(p_{\nu +r-2})\ll p_\nu \quad
\forall\,j .
\label{v2}
\end{equation}
(We may note that
$N_{\nu +r-1} = \max_{j} |n_{j,\nu +r-1} |
\asymp p_\nu $.)

Now from (\ref{v1}),  (\ref{v2}) we have
$$
   0<\zeta (n_\nu)\ll\psi ( \kappa N_{\nu +r-1}).
$$
We may suppose
$\zeta (n_{\nu+1})<\zeta (n_\nu).$

The proof is complete.
\vskip+0.3cm

\par Theorem 1.4 follows immediately from Lemmas 1.1 and 1.2
as for the numbers constructed in Lemma 1.1 by  Lemma 1.2 we have approximations satisfying
(\ref{n3xaxa}) and in this case the inequality (\ref{n3xaxa}) is also valid for the best
apptoximations.
 \par Theorem 1.4 is proved.
 \vskip+0.3cm

\par
{\bf 1.5.\hspace{2mm} Proof of Theorem 1.3.} \vskip+0.3cm

\par
We need the following notation. Let $ r \ge 2 ,$ $ {\bf R}^{r+1}  $ be Euclidean  space
with Cartesian coordinates $ ( x_0, \dots , x_{r} )$ , \,\,$ {\bf Z}^{r+1} \subset {\bf
R}^{r+1}  $ be the lattice of integers,
 $ L^{r}  $ be the $r$-dimensional subspace in
$ {\bf R}^{r+1}$ orthogonal to the vector $ ( 1, \alpha _1, \dots , \alpha _r ) $ and the
$r$-tuple
 $ \alpha _1,\dots
,\alpha _r  $ satisfies (\ref{n3}) with
\begin{equation}
 \psi (y) =
e^{-\gamma y } ,\hspace{2mm}  \gamma  \in (0;1) .
        \label{z1}
\end{equation}
(So for our proof we need only ordinary Hinchin's singular linear forms rather than the
generalization from Theorem 1.4.)

Let $ {\bf R}^{r+2} = {\bf R}^{r+1} (x_0, \dots , x_{r} ) \times {\bf R}^1(z)  $ be the
product of
 $ {\bf R}^{r+1} $ and $ {\bf R}^1 $,
$$
 {\cal L }^{r+1} = L^{r} \times {\bf R}^1 ,
$$
$$
 {\cal L }^{r+1}_\delta  = \{ X \in {\bf R}^{r+2} :\hspace{2mm}
\rho (X, {\cal L}^{r+1} ) \le \delta  \},
$$
$$
E^{r+2}_H =\{ X = ( x_0,\dots , x_r, z ) : \hspace{2mm}
 \max \{ |x_0|, \dots , |x_r|, |z| \} \le H \},
$$
$$
\Pi ( \sigma ; H ) = E^{r+2}_H \cap {\cal L}^{r+1}_\delta ,
$$
$$
{\cal K } = \bigcup_{ t \ge 1 } \hspace{2mm} \Pi ( 2 e^{ - (\gamma -\varepsilon )t} ; t ),
\hspace{2mm} \varepsilon \in (0,\gamma ).
$$
The infinite domain $ {\cal K } $ has finite volume:
$$
Vol \hspace{2mm}  {\cal K } \ll \int_1^\infty\hspace{2mm}  t^{r+1} e^{-(\gamma
-\varepsilon )t} dt  < + \infty.
$$
Moreover from our choice of $\psi $ by
 (\ref{z1})
we have $ m_\nu  \in {\cal K } \hspace{2mm}  \forall\nu  \ge \nu _0 .$
\par Let $ B_\epsilon  $  be  a $(r\! +\! 2)$-dimensional
ball with radius
 $ \epsilon <
1/2 $ centered at
  $ ( 0, \dots , 0 , 1 ) \in {\bf R}^{r+2}. $
For a point
 $\xi
\in B_\epsilon $ we put in correspondence the
 $ (r\! +\! 1)$-dimensional lattice $ \Lambda _\xi = {\bf Z}^{r+1} \oplus \xi
{\bf Z} $ generated by
 $ {\bf Z}^{r+1} $ and
the point
 $ \xi . $
Let $T_\xi: {\bf R}^{r+2}\to {\bf R}^{r+2}$ be the
 linear transformation preserving
the lattice ${\bf Z}^{r+1}$ and transforming the vector $\xi$ into the unit vector
 $(0,...,0,1)$.
Consider the $(r\! +\! 1)$-dimensional subspace $T_\xi {\cal L}^{r+1}$ and define $\alpha
_{r+1}$ in such a way that the vector $(1,\alpha _1,...,\alpha _r,\alpha _{r+1})$ is
orthogonal to the subspace $T_\xi {\cal L}^{r+1}$.

\par
The proof of  the following lemma is is in general  similar to the original proof of the
Minkowski-Hlawka theorem (see for example \cite{Gruber}). It is based on well-known metric
procedure.

\vskip+0.3cm
\par{\bf Lemma 1.3.} \hspace{2mm}
{\it For almost all points
 $ \xi \in B_\epsilon $
(in  sense of Lebesgue measure)
 we have
\par 1) \hspace{2mm} $\Lambda _\xi \cap {\cal L}^{r+1} = \{ 0 \}, $
\par 2) \hspace{2mm}  the intersection $ \Lambda_\xi \cap {\cal K} $
contains the points
 $ m_\nu  $
 and at most finite number of other integer  points. } \vskip+0.3cm
\par {\bf Corollary.} \hspace{2mm}
{\it Almost all but a finite number of the best approximations for the $(r\!+\!1)$-tuple
$(\alpha _1,...,\alpha _r, \alpha _{r+1} )$ from the lattice ${\bf Z}^{r+2}$ coincide with
the best approximations for the $r$-tuple $(\alpha _1,...,\alpha _r ) $ from the lattice
${\bf Z}^{r+1}$.} \vskip+0.3cm
\par Proof of Lemma 1.3.

Let $ k  $ be natural, $ e_j  $ be  unit vectors in $ {\bf R}^{r+1} $, $ \chi(X) $ be the
characteristic function of the domain $ {\cal K }. $ We consider the value
 $$
S_\xi ( T ) = \sum_{k =1}^T \hspace{2mm} \sum_{ m} \chi ( m_0 e_0 + \dots + m_{r} e_{r} +
m_{r+1} \xi )   \ge 0 , $$ where the inner sum is  taken over all
$$ (m_0, \dots , m_{r  }  )\in {\bf Z}^{r+1}, \hspace{2mm} m_{r+1} \in {\bf Z}
\setminus \{ 0\}:  \hspace{2mm} \max \{ |m_0|, \dots , |m_{r+1}| \} = k .
$$
This sum calculates the number of points of the lattice $ \Lambda _\xi  $ with norm not
greater than  $ T$  lying in $ {\cal K} $ and different from the points of $ {\bf Z}^{r+1}
\subset {\bf R}^{r+1}$:
$$ S_\xi ( T ) = \# \{  m
=  m_0 e_0 + \dots + m_{r} e_{r} + m_{r+1} \xi,
$$
$$
 m_{r+1} \neq 0 ,\,\,\,
 0 < \max \{ |m_0|, \dots , |m_{r+1}| \} \le T \}.
$$
Observe that
$$ \int_{B_\epsilon} S_\xi ( T ) d\xi = \sum_{k =1}^T \hspace{2mm}
\sum_{ m }  \hspace{2mm} Vol ( B_\epsilon( m ) \cap {\cal K} ),
$$
where
$$
B_\epsilon( m ) = \{ X  = m_0 e_0 + \dots + m_{r} e_{r} + m_{r+1} \xi : \hspace{2mm} \xi
\in B_\epsilon \}
$$
It is clear that for
 $ \max|m_j| = k $
we have
$$
Vol (B_\epsilon( m ) \cap {\cal K} ) < Vol ( {\cal K} \cap \{  z \in R^{r+2}:\hspace{2mm}
\max_j |z_j | \ge  k/2 \} \ll e^{-\gamma_1  k },
$$
where $  0 < \gamma _1 < \gamma $. Hence for any
 $T$
$$
\int_{B_\epsilon} S_\xi ( T ) d\xi \ll \sum_{k =1}^T k^{r+2}   e^{-\gamma_1  k
 } \ll 1.  $$
Now we use  Levi's theorem to establish that for almost all
 $ \xi \in B_\epsilon $
there exists a finite limit $ \lim_{T \to \infty}
 S_\xi (T) $.
This means that for almost all
 $ \xi $
the intersection
  $ \Lambda_\xi \cap
{\cal K} $ consist of at most finite number of points different from
 $ m_\nu .$ \par
The proof is complete. \vskip+0.3cm

\par Now Theorem 1.3
can be proved by induction. For  $ r = 2 $ we have Hinchin's singular vector $ ( \alpha
_1, \alpha _2) $ satisfying the singularity condition with $ \psi (y ) = e^{-y}. $ The
induction step is performed in Lemma 1.3 and the proof is complete.

\vskip+0.7cm
 \centerline {\S 2.\hspace{2mm} \bf The best simultaneous Diophantine approximations.}

\vskip+0.3cm
  {\bf2.1. \hspace{2mm} Definitions.}
\vskip+0.3cm

For a $s$-tuple of real numbers $ \alpha = ( \alpha _1, ... , \alpha _s ) \in {\bf R}^s$
we define {\it the best simultaneous approximation} (briefly, b.s.a) as an integer point $
\zeta  = ( p,a_1, ... , a_s ) \in {\bf Z}^{s+1} $ such that
$$
D(\zeta ): = \max_{j=1,...,s} |p \alpha _j - a_j | < {\rm min}^* \max_{j=1,...,s} |q
\alpha _j - b_j |,
$$
where $ \min^* $ is taken over all $q,b_1,...,b_s$ under conditions
$$
1 \le q \le p   ; \hspace{2mm}
 (b_1, ... , b_s) \in {\bf Z}^s \setminus \{(a_1, ... , a_s )\}.
$$
In the case
 $ \alpha _j \not \in {\bf Q}$
all b.s.a. to
 $ \alpha $
form  infinite sequences
$$
 \zeta^\nu   = ( p^\nu ,a_1^\nu , ... , a_s^\nu  ) ,\hspace{2mm}
 \nu = 1,2, ...,
$$
$$ p^1 <...< p^\nu < p^{\nu +1} < ... $$
and
 $$
 D(\zeta^1 ) > ... > D(\zeta^\nu  ) > D(\zeta^{\nu +1} )  > ...  .
 $$
 Let
$$
 M_\nu [ \alpha ]  =
 \left(
\begin{array} {cccc}
p^\nu  & a_1^\nu  &  \dots & a_s^\nu  \cr \dots & \dots & \dots & \dots \cr p^{\nu +s} &
a_1^{\nu +s} & \dots & a_s^{\nu +s}
\end{array}
  \right)
$$
and  $ {\rm rk}\,\, M_\nu [\alpha ]   $ be the rank of the matrix $ M_\nu [\alpha ]   $.
Natural  number $ R(\alpha ) $, $2 \le R(\alpha ) \le s+1$ is defined as follows
$$
R(\alpha ) = \min\hspace{2mm}  \{ n  :\hspace{2mm}{\it there} \hspace{2mm} {\it
exists}\hspace{2mm}{\it a}\,\,\,{\it  lattice}\hspace{2mm}  \Lambda \subseteq {\bf
Z}^{s+1}, \,\, {\rm dim } \,\,\Lambda  = n
$$
$$
\hskip+3.0cm \hspace{2mm} \hspace{2mm} \hspace{2mm} \hspace{2mm} {\it and\hspace{2mm}  }
\nu _0 \in {\bf N} {\it \hspace{2mm} such \hspace{2mm} that}\hspace{2mm} for \hspace{2mm}
all\hspace{2mm} \nu  > \nu _0\hspace{2mm} \zeta ^\nu  \in \Lambda   \hspace{2mm} \} $$

The value $ {\rm dim}_{\bf Z} \,\,\alpha  $ is defined as the maximum number of reals $
\alpha _{i_1}, ... , \alpha _{i_m} $  chosen from $ ( \alpha_0 =1,  \alpha _1, ... , \alpha
_s ) \in {\bf R}^{s+1}$ to be linearly independent over
 ${\bf Z}$.

\vskip+0.3cm {\bf Proposition 1.}  \hspace{2mm} {\it For  $ s= 1$ and  any $ \nu  $ we have
the equality $ {\rm det}\,\,M_\nu [\alpha ] = \pm 1 $ (it implies that for any $  \nu   $
we have ${\rm rk}\,\, M_\nu [\alpha ] = 2 $).}

\vskip+0.3cm

{\bf Proposition 2.}  \hspace{2mm} {\it For any $s \ge 1$  the following equality is valid
$$
R(\alpha ) = {\rm dim}_{\bf Z}\,\, \alpha  .
$$}

\vskip+0.3cm

{\bf Proposition 3.}    \hspace{2mm} {\it Let $ s = 2 $ and $ \alpha _1, \alpha _2 $
together with $1$ be linearly independent over ${\bf Z}$. Then there exist infinitely many
naturals
 $ \nu $
such that $ {\rm rk} M_\nu [\alpha ] = 3 = {\rm dim}_Z \alpha $ (hence
 the inequality
$ {\rm det}\,\,M_\nu [\alpha ] \neq 0 $ holds for infinitely many values of $\nu$.).}

Propositions 1- 3 are well-known and can be easily verified (compare \cite{Lag2}).

\vskip+0.3cm
  {\bf2.2. \hspace{2mm} Counterexample to Lagarias' conjecture.}
\vskip+0.3cm

We  formulate our result from \cite{MLAG} which deals with
 the degeneracy of the dimension of the spaces
generated by  successive  b.s.a. It gives a counterexample to Lagarias' conjecture
\cite{Lag2}. We would like to point out that this result was obtained due to discussion
with Nikolai Dolbilin. We shall give a sketched proof  in  next two sections.

\vskip+0.3cm

{\bf Theorem 2.1}\hspace{2mm}    {\it Let $ s \ge 3 $. Then there exists an uncountable set
of $s$-tuples $ \alpha $ with  components $ \alpha _1, ... , \alpha _s $ linearly
independent together with $1$ over $\bf Z$ such that
$ {\rm  rk}\,\, M_\nu [\alpha ]  \le
3 \hspace{2mm} \forall \, \nu  \in N. $
(Hence for all  $\nu $ the equality $ {\rm det}
M_\nu [\alpha ] = 0$ is valid.)}

\vskip+0.3cm
 2.3.\,\,\, { \bf Inductive lemma for Theorem 2.1.}
\vskip+0.3cm We consider Euclidean space ${\bf R}^{s+1}$ with Cartesian coordinates
  $(x,y_1,...,y_s)$.
Letter $\ell$ will denote a ray from the origin of coordinates located in the half-space $
\{x >0\}$. For such ray  $\ell$ and for small enough positive $\epsilon$ the opened cone $
K_\epsilon (\ell) $  consists of all rays $\ell '$ such that the angle between
  $\ell$ and $\ell '$ is less than $\epsilon$.
For a point $\xi$ from the half-space
 $\{ x>0\}$
we define   $\ell (\xi) $ to be the ray $\{ \kappa \xi: \kappa \ge 0 \} $. Subspace $ \pi
\subseteq {\bf R}^{s+1} $ is defined to be {\it  absolutely rational
 }
if the lattice
  $ \Lambda = \pi \cap {\bf
Z}^{s+1} $ has  dimension equal to the dimension of the whole $\pi $: $ {\rm dim}\,\,
\Lambda = {\rm dim} \,\, \pi $. Let $\ell$ be a ray parallel to a vector
$(1,\beta_1,...,\beta_s)$. The best approximation to the ray $ \ell$  is defined as a point
 $\zeta \in {\bf Z}^{s+1}$
which is the b.s.a. to $\beta$.

In the case when each  $\beta_j $ is not a half of an integer  the sequence of all best
approximations
$$
\zeta^\nu   = (p^\nu ,a_1^\nu  ,..., a_s^\nu )  ,\,\,\,  p^1 < ... ,
p^\nu <p^{\nu +1} <...
$$
 to the ray $\ell$ is defined correctly. It is finite in the case
when  there exists an integer point on the ray $\ell$
 different from the origin and is infinite
in the opposite case. This sequence of the best approximations we  write as
$$
{\cal B} (\ell ) = \{\zeta ^1,\zeta ^2, ... , \zeta ^\nu , ...\}.
$$
Moreover, we use the following notation:
$$ {\cal B}_k^t (\ell ) =
\{\zeta ^k,\zeta ^{k+1}, ... , \zeta ^t\}.
$$

\vskip+0.3cm

{\bf Lemma 2.1.}\hspace{2mm} {\it Let $\Lambda   ={\bf Z}^{s+1}\cap \pi $ be a lattice
located in an absolutely rational subspace $\pi $, ${\rm dim}\,\,\pi \ge 2$. Let a point
$\zeta \in \Lambda  $ satisfy the condition
$$
{\cal B} (\ell (\zeta ) ) = \{\zeta ^1,\zeta ^2, ... , \zeta ^\tau  , ...,\zeta^t
\},\hspace{2mm} \zeta^t = \zeta
$$
and in addition
$$ {\cal B}_\tau ^t (\ell (\zeta ) ) =
\{ \zeta^\tau  , ...,\zeta^t \} \subset \Lambda \subset \pi .
$$
Let $D(\zeta _\tau ) < D(\xi )$ for any integer point
 $\xi $
which does not belong to
 $\pi $.

Then for some  $\epsilon > 0 $ any ray      $ \ell ' \subset K_\epsilon (\ell (\zeta ))$
satisfy the following conditions:

1) \hspace{2mm} $ {\cal B}(\ell ') \supset
 {\cal B}(\ell (\zeta ) )$;

2)\hspace{2mm} the sequence of  the best approximations to the ray $\ell '$  between the
approximations $\zeta^\tau  $ ¨ $\zeta^t$   lies completely in the subspace $ \pi $. }

We must remember that the all points in
$$
{\cal B}_\tau ^t (\ell (\zeta ) )
$$
obviously belong to the considered  sequence of  the best approximations to the ray $\ell
'$ between the approximations $\zeta^\tau  $ and $\zeta^t$ but it may happen that a number
of
  new
points appear.

Lemma 2.1 follows from two easy observations:

1.\hspace{2mm} for a small perturbation $\ell '$ of the ray $\ell$ the first best
approximations to  $\ell$ remains to be the best approximations to $\ell '$;

2. \hspace{2mm} a small  perturbation of $\ell$ does not enable   integer points not
belonging to    $\pi $ to become best approximations between the approximations
$\zeta^\tau $ and $\zeta^t$

\vskip+0.3cm
 {\bf 2.4. Sketch of the proof of Theorem 2.1.}
\vskip+0.3cm The proof uses two inductive steps.

The first step.\hspace{2mm} Applying Lemma 2.1 many times we construct absolutely rational
subspaces
$$
\pi _1,\rho _1,\pi _2,\rho _2,..., \pi _s,\rho _s
$$
with dimensions $ { \rm dim}\,\,\, \pi _j = 2$, $ { \rm dim}\,\,\, \rho _j = 3$ and a ray
$\ell = \ell (\zeta )$, \hspace{2mm} $ \zeta  \in \pi _s $ such that

A)\hspace{2mm} $ \pi _j, \pi _{j+1} \subset \rho _j, $

B)\,\,\, $ {\cal B}(\ell) = \{ \zeta^1,...,\zeta^{\tau _1}, \zeta^{\tau
_1+1},...,\zeta^{t_1}, \zeta^{t_1+1},...,\zeta^{\tau _2},..., \zeta^{\tau
_s+1},...,\zeta^{t_s}\},
$

where
$$\zeta^{t_s}=\zeta , \,\,\, t_0 =1,\,\,\, t_j-\tau _j \ge s+1, \hspace{2mm} \tau _j -
t_{j-1} \ge s+1 \hspace{2mm} \hspace{2mm} \forall \hspace{2mm} j
$$
and
$$ \zeta^{\tau_j+1},...,
\zeta^{t_j} \in \pi _j \hspace{2mm} \hspace{2mm} \forall \hspace{2mm} j , \,\,
\zeta^{t_j+1},...,\zeta^{\tau _{j+1}} \in \rho _j \hspace{2mm} \hspace{2mm} \forall
\hspace{2mm} j, $$

C) ${\cal B}(\ell)$ (as well as the union $ \bigcup_{j=1}^s \rho _j $ ) does not belong to
any
 $s$-dimensional subspace of $ {\bf R}^{s+1}$.

We perform the construction of such a ray $\ell (\zeta )$ for which the best
approximations admit A), B), C)  by means of Lemma 2.1 by an inductive procedure.

The beginning of the inductive procedure is trivial. Let the subspaces
$$
\pi _1,\rho _1,\pi _2,\rho _2,..., \pi _k,\rho _k
$$
and the ray $\ell (\zeta^{t _k}), \zeta^{t _k} \in \pi _k$ be already constructed. Then by
Lemma 2.1 we take $\zeta ^{\tau_{k+1}}$ with required  properties and choose absolutely
rational subspace
 $\pi _{k+1}$
such that the ray
 $\ell (\zeta ^{\tau _{k+1}}) $
lie in this subspace and the dimension of the subspace generated by all subspaces $ \pi
_1,\rho _1,..., \pi _{k+1} $ is maximal. Then by Lemma 1.4 in   $\pi _{k+1}$
we find a
point $\zeta^{t_{k+1}}$ with the requested properties.

The second step.\hspace{2mm} We must apply the procedure of the first step many times and
construct a sequence of rays
 $ \ell^k= \ell (\xi ^k),$
$ \xi ^k \in {\bf Z}^{s+1}$,  $k = 1,2,...$ in such a way that for any ray
 $ \ell^k$
the set of the best approximations ${\cal B}(\ell^k)$ consists of  $k$  successive blocks.
Each of these blocks must satisfy the conditions A), B), C) from the first step.

The limit ray for the sequence of rays
 $\ell^k$ will correspond to the numbers
$\alpha _1,...,\alpha _s$ with the requested in Theorem 2.1 properties: all successive $(s
+ 1)$ b.s.a. for $\alpha _1,...,\alpha _s$ will lie in two- or three-dimensional subspaces
and proposition 2 form Section 2.1 and the property C) lead to the independence of the
reals $1, \alpha _1,...,\alpha _s$ over rationals.

\vskip+0.3cm

{  \bf 2.5.\hspace{2mm} Best simultaneous approximations in different norms.}

\vskip+0.3cm

We consider a convex O-symmetric star function $f : {\bf R}^n \to {\bf R}_+$ satisfying
the conditions

1) $f$ is continuous,

2) $ f(x) \ge 0 \, \forall x \in {\bf R}^n,\hspace{2mm} f(x) = 0 \Longleftrightarrow x = 0
$,

3) $ f(-x) = f(x) \,\forall x \in {\bf R}^n$,

4) $ f(tx) = t f(x)\hspace{2mm} \forall x \in {\bf R}^n,\hspace{2mm} \forall t \in {\bf
R}_+$,

5) the set $ B_f^1 = \{ y \in {\bf R}^n :\hspace{2mm} f(y) \le 1\} $ is convex  and $ 0
\in {\rm int } B_f^1$.

It is well known (see \cite{CASGEOMqw}) that $f$ determines a norm in  ${\bf R}^n$.
Function (or norm)  $f$ is
 {\it strictly convex}  if the set $
B_f^1$ is strictly convex; that is, the boundary
 $ \partial B_f^1$
do not have segments of straight lines. We  use $B^\lambda _f(a) $ for the set
$$ B^\lambda _f(a) =
\{ y \in {\bf R}^n : \hspace{2mm} f(y-a) \le \lambda \},
$$ so $ B_f^1 = B_f^1(0)$.

For an $n$-tuple
 $ \alpha  = (\alpha _1,..., \alpha _n) \in {\bf R}^n$
we define
 {\it $f$-best simultaneous approximation}
($f$-b.s.a.) as an integer point $\tau  = (p, a_1,..., a_n) \in {\bf Z}^{n+1}$ such that
$p \ge 1$  and
$$
f(\alpha q - b ) > f(\alpha p - a)
$$
for all
$$
(q, b_1,..., b_n) \in {\bf Z}^{n+1}, \hspace{2mm} \hspace{2mm} 1\le q \le p-1
$$
and for all
$$
(p, b_1,..., b_n) \in {\bf Z}^{n+1} , \hspace{2mm} \hspace{2mm} b \neq a.
$$
In the case when
 $f$  determines the cube
$$ B_f^1 = \{ y = (y_1,...,y_n)\in {\bf R}^n: \max_j |y_j| \le 1\}
$$
our  definition leads to the classical definition of the b.s.a. considered in  previous
sections.

All the $f$-b.s.a. for  $\alpha $ form the sequences
$$
\tau _\nu  = (p_\nu , a_\nu   )\in { \bf Z}^{n+1}, \hspace{2mm} p_\nu  \in {\rm
N},\hspace{2mm} a_\nu  = ( a_{1,\nu },..., a_{n,\nu }) \in {\bf Z}^n,
$$
$$  p_1 < p_2 <... < p_\nu  < ...,
$$
$$
f(\alpha p_1 - a_1 ) > f(\alpha p_2 - a_2) > ... . f(\alpha p_\nu  - a_\nu  ) > ...,
$$
and these sequences are finite in the case
 $ \alpha  \in {\bf Q}^n$ and are infinite in the opposite situation.

Let $\xi_\nu  = (\xi_{1,\nu },..., \xi_{n,\nu })$ denotes the remainder vector $
\xi_{j,\nu } = \alpha_j p_\nu  - a_{j,\nu } $. Let
$$
\Xi_\nu = (\Xi_{1,\nu },..., \Xi_{n,\nu });\hspace{2mm} \Xi_{j,\nu } = \xi_{j,\nu } /
f(\xi_\nu );
$$
obviously,
 $ \Xi_\nu \in B_f^1$.
For a given vector
 $ \xi \in {\bf R}^n$
we also  use the notation $ \Xi (\xi ) = \xi / f (\xi ) \in B_f^1$. Moreover, for the
integer vector $ \zeta  = (p, a_1, ... , a_n) \in {\bf R}^{n+1}$ we  use the notation $
\xi^\alpha  (\zeta )   = (p\alpha _1 - a_1,..., p\alpha _n - a_n) \in {\bf R}^n$.

\vskip+0.3cm

 {\bf 2.6.\hspace{2mm}
The order of the best approximations.} \vskip+0.3cm

From the Minkowski convex body theorem applied to the  cylinder
\begin{equation}
\Omega_\nu   = \{ z = (x, y_1,...,y_n) \in {\rm R}^{n+1}:\hspace{2mm} |x| < p_{\nu
+1},\hspace{2mm} f(\alpha x - y) < f(\xi_\nu )\} \label{rwe1}
\end{equation}
(it does not contain nontrivial integer points) it follows that for any  $\nu $ one has
\begin{equation}
f(\xi_\nu ) \le C_1(f) p_{\nu +1} ^{-1/n} \label{3qw}
\end{equation}
with constant $ C_1(f) = 2 /({\rm Vol } B_f^1)^{1/n} $. On the other hand, we can show that
the following result is valid.

\vskip+0.3cm {\bf Theorem 2.2.}\hspace{2mm} {\it Let ${\rm dim}_{\bf Z} (1, \alpha _1,
..., \alpha _n) \ge 3$. Then
\begin{equation}
f(\xi_\nu )  p_{\nu +1} \to +\infty,\hspace{2mm} \nu  \to  +\infty . \label{u1qw}
\end{equation}
}

Proof.

1) Let $ \Lambda^2 \in {\bf Z}^{s+1} $ be a two-dimensional sublattice   and $ {\rm det}_2
\Lambda^2$ be the area of its fundamental domain. The set of all sublattices
$$
\{ \Lambda^2 \subset {\rm Z}^{s+1}:\hspace{2mm} {\rm det}_2 \Lambda^2 \le \gamma \}
$$
is finite for any  $\gamma $.

2) Consider a two-dimensional lattice $ \Lambda^2_\nu  = \langle \tau _\nu , \tau _{\nu +1}
\rangle_{\rm Z} $. From $ {\rm conv} (0,\tau _\nu , \tau _{\nu +1}) \subset \Omega_\nu $
it follows that
$$
\frac{1}{2} {\rm det}_2 \Lambda^2_\nu = {\rm vol}_2 ({\rm  conv }(0,\tau _\nu , \tau _{\nu
+1} )) \ll f(\xi_\nu ) p_{\nu +1}.
$$

3) From ${\rm dim}_{\bf Z} (1, \alpha _1, ..., \alpha _s) \ge 3$ it is easy to  deduce (see
proposition 2 from Section 2.1) that the sequence of all $f$-b.s.a. cannot asymptotically
lie in  a two-dimensional sublattice and  hence
for a fixed sequence of naturals
$\nu_k$ the embedding
$\cup_{k  = 1}^\infty \tau _{\nu_k} \subset
\cup_{n =1}^{  \nu _0} \Lambda^2_\nu $ never holds.

Theorem 2.2 immediately follows from 1), 2), 3).

We would like to refer to Hinchin once again as in \cite{HINS1}, \cite{HINS} he actually
proved that it is not possible to establish any specific rate of growth of the value
 $ f(\xi_\nu ) p_{\nu +1}$
in (\ref{u1qw}):

\vskip+0.3cm {\bf Proposition 4.}\hspace{2mm} {\it For any function $\psi (y) \uparrow
+\infty$ increasing to infinity (as slow as one wishes) as $y \to \infty$ there exists an
$n$-tuple
$$ {\bf \alpha} \in {\bf R}^n,\,\,
{\rm dim}_{\bf Z} (1, \alpha _1, ..., \alpha _n) = n+1$$ such that
\begin{equation}
 f(\xi _\nu )  p_{\nu +1} = O(\psi (p_{\nu +1 })),\hspace{2mm}
\nu  \to  +\infty . \label{u2qw}
\end{equation}
}

Formula (\ref{u2qw}) shows that in the situation
 $ n \ge 2$
there exist vectors $\alpha $ for which  the lower estimate from (\ref{3qw},\ref{u1qw}) is
the exact one. Of course, in the case $ n =1$  for any  $\nu $ we have
$$
C_2(f) p_{\nu +1}^{-1} \le f(\xi_\nu ) \le C_1(f) p_{\nu +1}^{-1}
$$
(see \cite{Hin}).

\vskip+0.3cm

{\bf  2.7.\hspace{2mm} The directions of the successive  best approximations.} \vskip+0.3cm

{\bf Theorem 2.3.}\hspace{2mm} {\it  For any natural $\nu $ one has $ \Xi_{\nu +1} \not\in
{\rm int} B_f^1 (\Xi_\nu )$.}

Theorem 2.3 was actually proved by Rogers in \cite{ROD1qw} for signatures (see Section
2.11). It follows from the fact that in the cylinder (\ref{rwe1}) there is no nontrivial
integer points and
$$ \tau _{\nu +1} - \tau _\nu
= (p_{\nu +1} - p_\nu , a_{1,\nu +1} - a_{1,\nu },..., a_{n,\nu +1} - a_{n,\nu })
$$
does not belong to the cylinder $\Omega_\nu  $. Then one must notice that
 $ 0 < p_{\nu +1} - p_\nu < p_{\nu +1}$.
Hence
 $ \tau _{\nu +1} - \tau _\nu \not\in \Omega$
means that
\begin{equation}
\xi_{\nu +1} \not\in {\rm int} B_f^{f(\xi_\nu )} (\xi_\nu ) \label{2qw}.
\end{equation}
Now $ 0 \in \partial S_f^{f(\xi_\nu )}(\xi_\nu )$ and due to convexity we have
$$
\xi_{\nu +1} \frac{f(\xi_\nu )}{f(\xi_{\nu +1})} \not\in {\rm int} B_f^{f(\xi_\nu )}
(\xi_\nu )
$$
and this is exactly what is stated in the theorem.

We can notice that the  statement (\ref{2qw}) is a little bit more general than the
Theorem 2.3.

\vskip+0.3cm

{\bf 2.8.     \hspace{2mm} Strictly convex norms.}

\vskip+0.3cm

{\bf Theorem 2.4.}\hspace{2mm} {\it Let the norm $f$ be strictly convex. Then there exists
$ \delta = \delta (f) >0   $ such that for any vector $ \alpha  \not\in  {\bf Q}^n$ there
exist infinitely many values of $\nu $ for each of them
$$
\Xi_{\nu +1} \not\in B_f^{1+\delta }(\Xi_\nu ).
$$
}

We remind the reader that the $n$-tuple
 $\alpha = (\alpha _1,..., \alpha _n)$
is defined to be
 {\it badly approximable}
if for some positive $ D(\alpha ) > 0$ the inequality
$$
\max_{1\le j \le n} \hspace{2mm} \min_{a_j \in {\rm Z}}\hspace{2mm} |p\alpha _j - a_j |
\ge D(\alpha )p^{-1/n}    \hspace{2mm}
$$
is valid for all $ p\in {\rm N} $ (Concerning  the existence of the badly approximable
vectors see \cite{SCHMIDT}). It is easy to see that the  vector  $\alpha $ is badly
approximable if and only if for any norm
 $f$
there is a constant $ D_1(f,\alpha ) $ such that for any natural $p$ holds
\begin{equation}
\min_{a\in {\bf Z}^n} f(p\alpha -a) \ge D_1(f,\alpha ) p^{-1/n}. \label{Bqw}
\end{equation}

{\bf Theorem 2.5.}\hspace{2mm} {\it Let $\alpha $ be badly approximable and $D=D(\alpha )$
be the corresponding constant. Let the norm $f$ be strictly convex. Then there exist $w =
w(D,f) \in {\bf N}$ and $ \delta = \delta (D,f) > 0$ with the following property:

for any $ \nu \ge 1\hspace{2mm} $ there exists a natural $j$ from the interval $ \nu  \le j
\le \nu +w$ such that
\begin{equation}
\Xi_{j +1} \not\in B_f^{1+\delta }(\Xi_j ). \label{Cqw}
\end{equation}
}

Theorem 2.4 shows that for a strictly convex norm the condition $\theta_{\nu +1} \not\in
{\rm int} B_f^1(\theta ) $ for the sequence  $ \theta_\nu  \in B_f^1$ is not sufficient
for the existence of
 $\alpha $  such that
  $ \lim_{\nu \to \infty} (\theta_\nu  - \Xi_\nu ) = 0$.
Theorem 2.5 shows that for the badly approximable numbers the values of $j$ for which we
have (\ref{Cqw}) appear regularly. Probably, the result of Theorem 2.5 does not depends on
the fact that $\alpha$ is badly approximable but we cannot prove it. The proofs of the
theorems we give in next two sections. From the results of the Section 2.11 it is clear
that Theorem 2.4 {\it is not valid   } for {\it non-strictly} convex norms. \vskip+0.3cm

{\bf 2.9. \hspace{2mm}  Two lemmas.} \vskip+0.3cm

Lagarias \cite{Lag1} proved the following statement. \vskip+0.3cm {\bf Lemma
2.2.}\hspace{2mm} {\it Define $ h = 2^{n+1}$. Then
 for any norm $f$ and for any  natural $\nu $ one has $ p_{\nu +h} \ge 2 p_\nu $.}

\vskip+0.3cm {\bf Corollary 1.}\hspace{2mm} {\it For any vector
 $\alpha \not\in {\bf Q}^n$
and for all
 $\nu , j >1$  one has
$$
f(\xi_{\nu +jh}) \le C_1p_\nu ^{-1/n} \left( \frac{1}{2}\right)^{j/n}.
$$
}

\vskip+0.3cm

{\bf Corollary 2.}\hspace{2mm} {\it Let $ \alpha $ be badly approximable. Then there exists
$ h^* = h^*(f,\alpha ) \in {\bf N}$ such that
\begin{equation}
\forall \, \nu  \ge 1 \hspace{2mm} f(\xi_{\nu +h^*}) < \frac{1}{2}\, f(\xi_\nu )
\label{6qw}
\end{equation}
} \vskip+0.3cm

Proof of the Corollary 2.

From (\ref{3qw}) and the  condition (\ref{Bqw}) it follows that
$$
D_1 p_\nu ^{-1/n} \le f(\xi_\nu ) \le C_1 p_\nu ^{-1/n}
$$

and by Lemma 2.2 $p_\nu $ grows exponentially. Now (\ref{6qw}) follows.

\vskip+0.3cm

{\bf Lemma 2.3. }\hspace{2mm} {\it Let $f$ be strictly convex. Then for any $\varepsilon
>0$ there exists $\delta >0$ such that for any $ \theta \in \partial B_f^1(0)$ and  any
$ \xi \in B_f^1(0)\setminus B_f^1(\theta )$ under condition
$$
 \Xi (\xi ) \in \partial B_f^1(0)
\bigcap \left( B_f^{1+\delta }(\theta ) \setminus B_f^1(\theta ) \right)
$$
we have $ f(\xi )  > 1 - \varepsilon$.} \vskip+0.3cm

Proof. \vskip+0.3cm

Let $ \eta \in \partial B_f^1 (\theta ) \bigcap \partial B_f^1 (0)$. As $f$ is strictly
convex we have $ (0;\eta ) \subset {\rm int } B_f^1 (\theta )$. Now if $\Xi \in \partial
B_f^1(0) \setminus B_f^1(\theta ) $ belongs to a small $\delta $-neighborhood of the point
$\eta$ then the segment $[0; \Xi ]$ must intersect with
 $\partial B_f^1 (\theta )$ in some point
$ \zeta (\Xi ) = [0;\Xi ] \cap ( \partial B_f^1(\theta ) \setminus 0)$ and $ \zeta (\Xi )
\to \eta $ when $ \Xi \to \eta$. If $\xi \in B_f^1(0) \setminus B_f^1 (\theta )$ then
 $\xi$ is between  $\Xi (\xi ) $ and
$ \zeta (\Xi (\xi ))$.

Lemma is proved.

\vskip+0.3cm

 {\bf 2.10.\hspace{2mm} The proofs of Theorems 2.4, 2.5.} \vskip+0.3cm

We  prove Theorem 2.4.

Suppose that Theorem 2.4 is not valid. Then for any $ \delta  > 0$ we have
$$
\Xi_{\nu +1} \in B_f^{1+\delta }(\Xi_\nu )
$$
 for
$ \nu \ge \nu _0 (\delta ). $ Now from Theorem 2.3 and Lemma 2.3 we deduce that for any $
\varepsilon > 0$
$$
f(\xi_{\nu +1} ) \ge ( 1- \varepsilon ) f(\xi_\nu  )
$$
when $ \nu  \ge \nu _0(\varepsilon ). $ It means that
\begin{equation}
f(\xi_{\nu _0 +j}) \ge (1-\varepsilon )^j f(\xi_{\nu _0}). \label{7qw}
\end{equation}
But from Corollary 1 to Lemma 2.2 we see that
\begin{equation}
f(\xi_{\nu _0 +j}) \le C_1 p_{\nu _0}^{-1/n} (1/2)^{j/n} . \label{8qw}
\end{equation}
For small values of
 $\varepsilon$  the inequalities (\ref{7qw}) and (\ref{8qw})
lead to contradiction when $ j \to \infty$.

Theorem 2.4 is proved.

\vskip+0.3cm

Now we  prove Theorem 2.5.

Suppose that Theorem 2.5 is not valid. In this situation for arbitrary large $ w \in {\bf
N}$ and for arbitrary small  $ \delta  > 0$ there exists $\nu $ satisfying the condition
$$
\Xi_{j+1}  \in B_f^{1+\delta } (\Xi_j),\hspace{2mm} j = \nu ,\nu +1,..., \nu +w.
$$
Applying Lemma 2.3 we see that
\begin{equation}
f(\xi_{\nu +w} ) \ge (1-\varepsilon)^w f(\xi_\nu ), \label{10qw}
\end{equation}
and  $\varepsilon >0$ may be taken arbitrary small. But at the same time from Corollary 2
to Lemma 2.2 we deduce that
\begin{equation}
f(\xi_{\nu +w} ) \le \left( \frac{1}{2} \right)^{[w/h^*]}
 f(\xi_\nu ).
\label{11qw}
\end{equation}
Again we take $\varepsilon$  small enough and the inequalities (\ref{10qw}) and
(\ref{11qw}) lead to contradiction when $ w \to \infty$.

The proofs are complete.

\vskip+0.3cm

{\bf 2.11. \hspace{2mm} Result on signatures and  illuminated points. } \vskip+0.3cm

For vector $\eta = (\eta_1,...,\eta_n)$ its signature is defined as
$$
{\rm sign }\,\eta = ( {\rm sign}\, \eta_1 ,..., {\rm sign}\, \eta_n ).
$$
Rogers \cite{ROD1qw} showed that for ordinary b.s.a. (in the case $ B_f^1 = \{ y: \max_j
|y_j| \le 1\}$) the successive best approximations satisfy the condition
$$
\forall \, \nu \hspace{2mm} {\rm sign }\, \xi_\nu  \neq {\rm sign }\,\xi_{\nu +1} .
$$
 (This
simple result was generalized in Theorem 2.3.) On the other hand, Sos and Szekeres
\cite{SOSqw} proved that for any sequence of signatures $\{ \sigma _\nu \}$ with
 $ \sigma
_\nu \neq \sigma _{\nu +1}$ there exists a vector $\alpha  \in {\bf R}^n$ with components
$\alpha _1,..., \alpha _n$ linearly independent together with $1$ over ${\bf Z}$  such
that $ {\rm sign}\, \xi_\nu = \sigma _\nu $. We give a generalization of this result.

Let $M\subset {\bf R}^n$ be a convex closed domain,
 $ b \in \partial M$
and $ a\not\in  M$. The point $b$ (as a point of the  boundary $\partial M$) is {\it
illuminated} from the point $a$  if there exists a positive $\lambda $ such that
$b+\lambda (b-a) \in {\rm int M}$.

\vskip+0.3cm

{\bf Theorem 2.6.}\hspace{2mm} {\it Let the sequence of points $\{ \theta_\nu \}_{\nu  =
1}^\infty \subset B_f^1$ satisfy the following condition: for each $\nu $ the point $0 $ as
the point of the boundary $ \partial B_f^1(\theta_\nu ) $ is illuminated from the point
$\theta_{\nu  +1}$. Then there exists  a vector $\alpha = (\alpha _1,..., \alpha _n)$ with
linearly independent components such that
\begin{equation}
\lim_{\nu \to +\infty} | \Xi_\nu  - \theta_\nu | = 0 . \label{limqw}
\end{equation}
} \vskip+0.3cm

{\bf Remark 1.}\hspace{2mm} {\it The result by Sos and Szekeres on signatures immediately
follows from our Theorem 2.6.}

{\bf Remark 2.}\hspace{2mm} {\it In (\ref{limqw}) we can provide any rate of convergence
to zero.}

We would like to say that the formulation of the conditions of Theorem 2.6 in terms of
illuminated points is due to
 O. German \cite{GER}. Moreover O. German
\cite{GER}
 proved
some interesting and new results on distribution of directions for
the best approximations in sense of linear forms and on the rate
of convergence to the asymptotic directions. In the next section
we shall give a sketch of the proof  of this theorem  and here we
consider one example.

In the case
 $ s = 2$
we consider the norm $ f^*(x_1, x_2)$ with unit ball $B_f^1$ defined by the inequalities
$$
| x_1 + x_2| \le 4, \hspace{2mm} |x_1 - x_2| \le 1.
$$
Applying Theorem 2.6 and observing the geometry of mutual configuration of the balls $
B_f^1 (0)$ and $  B_f^1 (\theta)$ we obtain the following statement. \vskip+0.3cm

{\bf Theorem 2.7.}\hspace{2mm} {\it For the norm $ f^*(x)$ the set of all  $f$-b.s.a. may
have the constant sequence of signatures: $ \sigma _\nu  = (+,+) , \hspace{2mm}  \forall
\nu $. } \vskip+0.3cm

Theorem 2.7 shows that the conclusion of Rogers' theorem from \cite {ROD1qw} is not true
for the norm $ f^*(x) $. One can easily construct the corresponding multi-dimensional
example and an example with strictly convex norm.

We would like to point out that we cannot construct an example of a norm $f$  for which
the sequence of all $f$-b.s.a. can have any given sequence of signatures. We may
conjecture that th Euclidean norm $ f(x_1, ..., x_n) = \sqrt{ x_1^2 +... + x_n^2} $ has
this property.

\vskip+0.3cm
 {\bf 2.12.\hspace{2mm}
Sketch of the proof of Theorem 2.6.  } \vskip+0.3cm

The proof is performed in the same manner as the proof of Theorem 2.1. By means of some
inductive procedure we construct a sequence of integer points
$$\tau _\nu
= ( p_\nu ,a_{1,\nu },..., a_{n,\nu })
$$
which must form the sequence of all $f$-b.s.a. for the limit point
$$
\lim_{\nu  \to +\infty} (a_{1,\nu }/p_\nu ,..., a_{n,\nu }/p_\nu ).
$$

The base of induction is trivial.

We sketch the induction step.

Let the points
$$
\tau _1,..., \tau _\nu  \in {\bf Z}^{n+1},\,\,\,
 \tau _j = (p_j; a_j) = (p_j,
a_{1,j},..., a_{n,j}),\,\,\,  1\le p_1 < p_2 < ... < p_\nu
$$
be constructed satisfying the following conditions:

1) \hspace{2mm} $ \tau _1,..., \tau _\nu  $ is the set of {\it  all} $f$-b.s.a. to
rational vector
$$\beta^\nu = (a_{1,\nu }/p_\nu , ..., a_{n,\nu } /p_\nu ),
$$

2)\hspace{2mm} $\Xi (\xi^{\beta^\nu }(\tau _j)) - \theta_j$ is small for all

$j = 1, ..., \nu - 1 $,

3)\hspace{2mm}
 $\theta_j
$ illuminates the point 0 of the boundary of $\partial  B_f^1 (\Xi (\xi^{\beta^\nu }(\tau
_{j-1}))$ for all $ j = 1,..., \nu $,

4) \hspace{2mm} there is no integer points on the boundary of the cylinder
$$
\{ (x,y_1,...,y_n) \in {\bf R}^{n+1}:\hspace{2mm} |x| < p_{\nu } ,\, f(\alpha x-y) \le
f(\xi_{\nu -1})  \}
$$
but the best approximations.

We must show how one can determine an integer point
$$ \tau _{\nu +1} = (p_{\nu +1}; a_{\nu
+1}) = (p_{\nu +1}, a_{1,\nu +1},..., a_{n,\nu +1}), \,\,\,  p_\nu  < p_{\nu +1}
$$
such
that

1$^*$) \hspace{2mm} $ \tau _1,..., \tau _{\nu +1} $ are {\it  all} $f$-b.s.a. to rational
vector
$$\beta^{\nu +1}= (a_{1,\nu +1 }/p_{\nu +1} , ..., a_{n,\nu +1 } /p_{\nu +1} ),$$

2$^*$)\hspace{2mm} $\Xi (\xi^{\beta^{\nu  +1}}(\tau _j)) - \theta_j$
 is small for all  $j = 1, ..., \nu $,

3$^*$)\hspace{2mm} $\theta_j $ illuminates the point 0 of the boundary $ B_f^1 (\Xi
(\xi^{\beta^{\nu +1} }(\tau _{j-1}))$ for all

$ j = 1,..., \nu +1 $,

4$^*$) \hspace{2mm} there is no integer points on the boundary of the cylinder
$$
\{ (x,y_1,...,y_n) \in {\bf R}^{n+1}:\hspace{2mm} |x| < p_{\nu +1} ,\, f(\alpha x-y) \le
f(\xi_\nu+1)  \}
$$
but the best approximations.

Consider a small neighborhood   of $ B^\lambda _f(a_\nu )$. Let $\lambda $ be small enough.
Then for any $\beta \in
 B^\lambda _f(a_\nu )$
integer points $ \tau _1,..., \tau _{\nu  -1} $  form {\it  all} the first successive
 $ \nu -1$ $f$-b.s.a to $\beta$. The  main difficulty is that for any
 $ \lambda $
there are some $\beta \in
 B^\lambda _f(a_\nu )$
for which between
 $ \tau _{\nu - 1} $  and $ \tau _{\nu }$
must arrive one new $f$-b.s.a. and it must be controlled.

We consider the ball $B^1_f$ and  the point $\Xi (\xi ^{\beta^\nu } (\tau _{\nu -1})) \in
\partial B^1_f$. Let $ B^* = B_f^1 (\Xi (\xi ^{\beta^\nu } (\tau _{\nu -1}))).$ From the
induction hypotheses 3) we know that $ \theta_\nu $ illuminates $0\in B^* $.

Then near the point
 $\theta_\nu \cdot t $  for some positive $t$  in the set
$ B_f^1 \cap {\rm int} B^*$
 there exists a point $\Xi^*$
such that the two-dimensional subspace $\pi ^*$ generated by the points $ \tau _\nu ,
\zeta ^* = (p_\nu , \Xi^*) $
 is
absolutely   rational. The point
 $\Xi^*$
must be very close to
 $ \theta_\nu \cdot t$.
So due to continuity we can obtain $\theta_{\nu +1}  \cdot t' \in {\rm int} B_f^1 (\Xi^*)$
for some positive $t'$.

In the absolutely rational subspace $\pi ^*$  a point $ \tau _{\nu +1} = (p_{\nu +1};
a_{\nu +1}) $ must be taken in such a way that $\tau _{\nu +1} $  and  $\zeta ^*$ be by
the same side of the line
 $0\tau _\nu $
(here we use the fact that $\pi ^*$  has dimension  2). As $\Xi^* \in {\rm int} B^*$ we
deduce from convexity that the whole segment $(0; \Xi^*)$ is in  ${\rm int}B^*$. Now we
can choose $\tau _{\nu +1}$ very close to the line  $0\tau _\nu $ and hence the sequence $
\tau _1,..., \tau _{\nu +1} $ really be the set of {\it  all} $f$-b.s.a. to the rational
vector $\beta^{\nu +1}= (a_{1,\nu +1 }/p_{\nu +1} , ..., a_{n,\nu +1 } /p_{\nu +1} )$.

The inductive step is sketched.

The sequence of vectors $\beta^\nu $ converges due to the smallness of the difference
$|\beta^{\nu +1} -\beta^\nu |$.

Linear independence over ${\bf Z}$ of the limit numbers
 $1,\alpha _1,..., \alpha _n$
may be obtained by the application of the Proposition 2 of Section 2.1.

\vskip+0.3cm

{\bf2.13.\hspace{2mm} Asymptotic directions.} \vskip+0.3cm

In this section we formulate without proofs some simple corollaries from our previous
results in terms of the asymptotic directions for the best approximations.

{\it The asymptotic direction  } for the
 $f$-b.s.a. sequence for a vector  $\alpha $
is defined as a point $ \theta \in \partial B_f^1(0)$ such that there exists  a subsequence
 $ \nu _j$
with the property $ \lim_{j \to +\infty} \Xi_{\nu _j} = \theta$. The set of all asymptotic
directions for $\alpha $ we  denote by $\Gamma _f(\alpha )$. Obviously
 $ \Gamma  _f(\alpha ) \subseteq B_f^1(0)$ is closed.

It seems to the author that C. Rogers was the first who gave the  definition of the
asymptotic direction for the Diophantine approximations \cite{ROD2qw} but our definition
differs from the Rogers' definition.

A set
 ${\cal A} \subseteq B_f^1(0)$
is defined to be {\it $f$-asymptotically admissible} if there is  an infinite sequence
$$
\theta_0 , \theta_1, ... , \theta_k ,...,\hspace{2mm} {\rm with }
 \hspace{2mm} \theta_j \in {\cal A}
$$
such that

1)  $\theta_j $ illuminates the point $0\in \partial  B_f^1 (\theta_{j-1})$,

2)   the set of all limiting  points of the sequence $\{\theta_k\}$ is just ${\cal A}$.

{\bf Theorem 2.8.}

  {\it Let  ${\cal A} \subseteq B_f^1(0)$
be  $f$-asymptotically admissible. Then there exists a vector $\alpha  \in {\bf R}^n$ with
linearly independent over rationals components such that ${\cal A} = \Gamma _f(\alpha )$.}

{\bf Corollary.}\hspace{2mm}    {\it If ${\cal A}$  is closed and there is $x \in {\cal
A}$ such that $-x \in {\cal A}$ then there exists
 $\alpha  \in {\bf R}^n$ with independent components such that
${\cal A} = \Gamma _f(\alpha )$.}

This result may be compared to Rogers observation \cite{ROD2qw} that the set of all
asymptotic directions is  not necessary 0-symmetric but there is some kind of symmetry.

The next theorem follows from Theorem 2.4.

{\bf Theorem 2.9.} \hspace{2mm} {\it Let norm $f$ be strictly
convex. Then there exists a positive $ \delta _1  $ depending on
 $f$
such that in the case
$$
{\cal A} \subset {\rm int} B_f^{1+\delta _1}(\theta )\hspace{2mm} \forall \, \theta \in
{\cal A}
$$
$ {\cal A}$ cannot be the set of the form $\Gamma_f (\alpha )$.}

\end{document}